\documentclass [12pt]{amsart}

\usepackage{amssymb,amsxtra,amsfonts}

\newcommand\beq{\begin{equation}}
\newcommand\eeq{\end{equation}}

\newcommand{\IP}{\mathbb{P}}                                     
\newcommand{\IQ}{\mathbb{Q}}                           
\newcommand{\IR}{\mathbb{R}}                           
\newcommand{\IC}{\mathbb{C}}

\newcommand{\cF}{\mathcal{F}}

\newcommand{\g}{       \mathfrak{g}     }

\newcommand{\PVI}{$\text{P}_{\text{VI}}$}   

\newcommand{\lsl}{\mathfrak{sl}} 

\newcommand{\pf}{\begin{bpf}}

\newcommand{\pfms}{\begin{bpfms}}
\newcommand{\epf}{\end{bpf}\hfill$\square$\\}           
\newcommand{\epfms}{\end{bpfms}\hfill$\square$\\}               

\newcommand{\idea}{\begin{bidea}}

\newcommand{\eidea}{\end{bidea}\hfill$\square$\\}           

\newcommand{\sk}{\begin{bsk}}    

\newcommand{\esk}{\end{bsk}\hfill$\square$\\}           
\newcommand{\sketch}{\begin{bsketch}}

\newcommand{\esketch}{\end{bsketch}\hfill$\square$\\}

\newcommand{\wt}{\widetilde}

\newcommand{\ga}{\gamma}

\newcommand{\Ga}{\Gamma}

\newcommand{\Sym}{\text{\rm Sym}}

\newcommand{\SL}{\text{\rm SL}}
\newcommand{\PSL}{\text{\rm PSL}}

\newcommand{\SO}{\text{\rm SO}}

\def\mapright#1{\smash{
        \mathop{\longrightarrow}\limits^{#1}}}

\def\mapdown#1{\Big\downarrow
        \rlap{$\vcenter{\hbox{$\scriptstyle#1$}}$}}

\theoremstyle{plain}
\newtheorem {hypo}{\bf\hspace{-\parindent}Hypothesis}
\newtheorem*{plainthm}{Theorem}

\newtheorem {prop}[hypo]{Proposition}

\newtheorem {cor}[hypo]{Corollary}

\theoremstyle{definition}
\newtheorem {defn}[hypo]{Definition}

\theoremstyle{remark}
\newtheorem {rmk}[hypo]{Remark}

\makeatletter
 \newlength{\typesize}
\makeatother

\newlength{\vvoff}
\newlength{\hhoff}

\newcommand{\locateoffcenter}[1]{%
\addtolength{\vvoff}{-0.25\typesize}%
\raisebox{\vvoff}{\hspace{\hhoff}\makebox(0,0){\smash{#1}}}
}
\newcommand{\object}[1]{%
\setlength{\vvoff}{0pt}%
\setlength{\hhoff}{0pt}%
\locateoffcenter{#1}
}

\newcommand{\swlabel}[1]{%
\setlength{\vvoff}{-0.5\typesize}%
\setlength{\hhoff}{0.75\typesize}%
\locateoffcenter{#1}
}
\newcommand{\nwlabel}[1]{%
\setlength{\vvoff}{-0.5\typesize}%
\setlength{\hhoff}{-0.75\typesize}%
\locateoffcenter{#1}
}

\begin{document}


\title[Higher genus Icosahedral Painlev\'e curves]
{Higher genus Icosahedral Painlev\'e curves}
\author{Philip  Boalch}
\address{\'Ecole Normale Sup\'erieure\\
45 rue d'Ulm\\
75005 Paris\\
France} 
\email{boalch@dma.ens.fr}


\begin{abstract}
We will write down the higher genus algebraic curves 
supporting icosahedral solutions of the sixth Painlev\'e equation,
including the largest (genus seven) curve.
\end{abstract}


\maketitle


\renewcommand{\baselinestretch}{1.02}            
\normalsize

\begin{section}{Introduction}

A Painlev\'e curve $\Pi$ is an algebraic curve supporting a solution to 
Painlev\'e's sixth equation (henceforth \PVI).
That is, there should be  rational functions
$y,t$ on $\Pi$ such that
\beq \label{eq: t cover}
t:\Pi\longrightarrow \IP^1\eeq
is a Belyi map (so expresses $\Pi$ as a branched cover, ramified 
only over $0,1,\infty$)
and $y$ (viewed as a function of $t$)
solves a \PVI\ equation.

This notion was introduced by Hitchin 
\cite{Hit-Poncelet} who found an infinite family of examples 
related to the Poncelet problem.
In essence he showed that all the modular curves $X_1(n)$ 
are Painlev\'e curves, at least 
for $n$ prime.
More precisely one should first pull back along the standard map
$X(2)\to X(1)$ (with Galois group $\Sym_3=\PSL_2(2)$),
 so there is a diagram:
$$
\begin{array}{cccc}
 & \Pi  & \mapright{} & X_1(n) \\
&\mapdown{t} &&  \mapdown{} \\
 \IP^1\cong\!\!\!\! &X(2)  & \mapright{} & X(1).
\end{array}
$$
In particular, for $n=5$, Hitchin wrote down the first 
explicit genus one Painlev\'e curve.

The aim of this article is to write down some other
explicit Painlev\'e curves not in
the above family of examples.

The (nonlinear) \PVI\ equation controls the ``isomonodromic'' 
(or monodromy preserving) deformations
of (linear) rank two Fuchsian systems on $\IP^1$ with four singularities, 
at $0,t,1,\infty$.
The monodromy of such a system is a representation
$$\rho:\cF_3=\pi_1(\IP^1\setminus\{0,t,1,\infty\})\ \mapright{}\ 
\SL_2(\IC)$$
and one of the main properties of these Painlev\'e curves is that 
the monodromy of the cover \eqref{eq: t cover},
i.e. its permutation representation $\cF_2\to \{1,2,\ldots,\deg(t)\}$,
coincides with the standard action of the pure mapping class group 
of the four-punctured sphere
($\cong \cF_2$) on the orbit it generates through the conjugacy class of the
representation $\rho$.

Hitchin's examples arose by seeking such isomonodromic 
deformations when the image 
of $\rho$  was equal to a binary dihedral group, and
in a previous article \cite{icosa} the author studied the case
when the monodromy group is equal to the binary icosahedral group.
All such solutions were classified and explicit formulae 
were written down for all but $8$ of the $52$ cases, including all those of
genus zero and most of the genus one cases.
(Five interesting cases had previously appeared in 
\cite{Dub95long,DubMaz00,Kitaev-dessins}.)

Unfortunately the icosahedral Painlev\'e curves of genus $\ge 2$
were not amenable to the method of construction used in \cite{icosa},
essentially due to the large degrees of the Belyi maps $t$.
(The method used was to first obtain, from the icosahedral linear monodromy, 
the precise asymptotics of the \PVI\ solution, using 
(the author's correction of)
Jimbo's asymptotic formula; this determined the Puiseux expansions to
arbitrary order which in turn enabled the curve to be obtained
algebraically.)

However it turns out that there is a trick to convert earlier icosahedral 
Painlev\'e curves  (that were found in \cite{icosa}, or were previously known)
into those of higher genus.
Namely one may use the so-called ``quadratic transformations''
introduced by Kitaev \cite{Kitaev-quad-p6} in 1991 and 
written in simpler form by Ramani et al. \cite{RGT-quad}
(we learnt of them from the recent article \cite{TOS-folding}).
Somewhat miraculously the solutions that can be obtained in
this way are almost exactly the complement of those 
we were able to obtain by the previous
method (there is a small overlap though).

Thus our aim
is to explain how the quadratic transformations may be applied in
this way and write down the resulting curves.
(This is not entirely trivial since, if applied blindly, the 
quadratic transformations lead to badly parameterised solutions, for example
with the wrong genus.) We also make some effort to obtain nice models (over
$\IQ$) of the resulting Painlev\'e curves.

For example the following result will be established:
\begin{plainthm}
There are precisely two non-hyperelliptic icosahedral Painlev\'e curves.
The first supports two inequivalent Painlev\'e solutions and 
is of genus three  and isomorphic to the smooth plane quartic with
affine equation
$$5(p^4+q^4)+6(p^2q^2+p^2+q^2)+1=0.$$
The second is of genus seven and is birationally isomorphic over $\IQ$ 
to
the affine curve cut out by the octic
$$
9\,(p^6\,q^2+p^2\,q^6)+
18\,p^4\,q^4+
4\,(p^6+q^6)+
26\,(p^4\,q^2+p^2\,q^4)+
8\,(p^4+q^4)+
57\,p^2\,q^2+
20\,(p^2+q^2)+
16
$$
whose closure in $\IP^2$ only has double point singularities. 
Moreover the obvious symmetries of these curves (negating and exchanging $p$
and $q$, generating a dihedral group of order $8$)
correspond to the Okamoto symmetries of the Painlev\'e solutions.
\end{plainthm}

\end{section}

\begin{section}{Background}

We will constrain ourselves to giving the notation 
and terminology that we will
use, referring the reader to \cite{icosa} or the review article 
\cite{srops} and references therein for more details and geometrical
background.

The sixth Painlev\'e equation (\PVI) is: 
\begin{align*}\frac{d^2y}{dt^2}=
&\frac{1}{2}\left(\frac{1}{y}+\frac{1}{y-1}+\frac{1}{y-t}\right)
\left(\frac{dy}{dt}\right)^2
-\left(\frac{1}{t}+\frac{1}{t-1}+\frac{1}{y-t}\right)\frac{dy}{dt}\\
&+\frac{y(y-1)(y-t)}{2\,t^2(t-1)^2}\left(
(\theta_4-1)^2-
\frac{\theta_1^2\, t}{y^2}+ 
\frac{\theta_3^2(t-1)}{(y-1)^2}+
\frac{(1-\theta_2^2)t(t-1)}{(y-t)^2}\right)
\end{align*}
where $\theta=(\theta_1,\theta_2,\theta_3,\theta_4)$ are (complex) constants.
This arises naturally when one tries to 
isomonodromically deform Fuchsian systems of the form
\beq\label{eq: lin syst}
\frac{d}{dz}-\left(
\frac{A_1}{z}+\frac{A_2}{z-t}+\frac{A_3}{z-1}\right),\qquad
A_i\in\g:=\lsl_2(\IC)
\eeq
as the second pole position $t$ varies in $\IP^1\setminus\{0,1,\infty\}$.
(The parameters $\theta$ specify  the eigenvalues of the residues: namely $A_i$
has eigenvalues $\pm\theta_i/2$ for $i=1,2,3,4$, where $A_4=-\sum_1^3 A_i$.)
Geometrically \PVI\ can (thus) be thought of as the explicit form of the
simplest nonabelian Gauss--Manin connection.

\begin{defn}
An algebraic solution of \PVI\ consists of a triple $(\Pi,y,t)$ where
$\Pi$ is a compact (possibly singular) algebraic curve and  $y,t$ are 
rational functions on $\Pi$ such that:

$\bullet$
$t:\Pi\to \IP^1$ is a Belyi map (i.e. $t$ expresses $\Pi$ as a branched cover
of $\IP^1$ which only ramifies over $0,1,\infty$), and

$\bullet$
Using $t$ as a local coordinate on $\Pi$ away from ramification points,
$y(t)$  should solve \PVI, for some value of the parameters $\theta$.
\end{defn}

Indeed given an algebraic solution in the form of a polynomial relation
$F(y,t)=0$ one may take 
$\Pi$ to be the closure in $\IP^2$ of the affine plane curve defined by $F$.
That $t$ is a Belyi map on $\Pi$ follows from the Painlev\'e property of \PVI:
solutions will only branch at $t=0,1,\infty$ and all other singularities
are just poles. The reason we prefer this reformulation is that often the
polynomial $F$ is quite complicated and usually there are much simpler
models of the plane curve defined by $F$. 
(The polynomial $F$ can of course 
be recovered as the minimal polynomial of $y$ over $\IC(t)$.)

We will say a Painlev\'e curve $\Pi$ 
is `minimal' or an `efficient parameterisation' if $y$ generates
the field of rational functions on $\Pi$, over $\IC(t)$, so that $y$ and $t$
are not pulled back from another curve covered by $\Pi$ (i.e. that
$\Pi$ is birational to the curve defined by $F$).

The main invariants of an algebraic solution are the genus of a (minimal) 
Painlev\'e
curve $\Pi$ and the
degree of the corresponding
Belyi map $t$ (the number of branches the solution has over the
$t$-line).

We will say that two solutions of \PVI\  are {\em equivalent}
if they are related by Okamoto's affine $F_4$ Weyl group symmetries
\cite{OkaPVI} of \PVI\  (which act on the set of parameters 
$\{\theta\}\cong\IC^4$ in the standard way).
(See e.g. \cite{TOS-folding,srops} for formulae for this action.)
For an algebraic solution, this acts within 
the set of rational functions on the
curve $\Pi$, and preserves the degree and genus of the solution
(at least if the linear monodromy representation is irreducible and not rigid).

We are interested here in the case where the monodromy group of 
the linear system \eqref{eq: lin syst} is equal to the binary icosahedral
group\footnote{more precisely we are interested in the solutions equivalent
  to such; one should bear in mind that the Okamoto transformations can change
  the monodromy group, and it will in fact be simpler to work at different
  equivalent values of the parameters $\theta$.
 cf. Remark \ref{rmk: maple file}}
$\Ga\subset \SL_2(\IC)$.
To understand the different cases that may occur 
essentially amounts to studying the different conjugacy classes of the local
{\em projective} monodromies.
Recall that the icosahedral rotation group $\Ga/\pm\cong A_5\subset\SO_3(\IR)$
has four non-trivial conjugacy classes, which we will label 
$a,b,c,d$ corresponding to rotations by 
$\frac{1}{2},\frac{1}{3},\frac{1}{5},\frac{2}{5}$-of a turn, 
respectively.
Thus we define, as in \cite{icosa}, the {\em $A_5$-type} of a representation
$$\rho:\pi_1(\IP^1\setminus\{0,t,1,\infty\})\to \Ga$$
to be the corresponding unordered 
set of four conjugacy classes of projective local monodromies 
(i.e. take the conjugacy classes of the images in $A_5$ of the
elements $\rho(\ga_i)$
for simple loops $\ga_i$ encircling one of $0,t,1$ or $\infty$ once).
The different cases that occur are tabulated in \cite{icosa}.

Two inequivalent 
icosahedral solutions will be said to be {\em siblings} if their 
monodromy representations $\rho$ are related by the nontrivial outer 
automorphism of $A_5$ (swapping the conjugacy classes $c,d$).
They will have the same Belyi map $t$, just a different solution function $y$.
(In general it is useful to
generalize this notion by considering Galois conjugate 
representations, e.g. for representations into  the $237$ triangle group 
there are sometimes three siblings, cf. \cite{octa}.)

\end{section}
\begin{section}{Quadratic transformations}

The basic idea \cite{Kitaev-quad-p6}
behind the quadratic transformations is as follows. 
Given an icosahedral Fuchsian system $A$ with $A_5$ type
$a^2\xi\eta$ for some $\xi,\eta\in\{a,b,c,d\}$
(i.e. with two local monodromies, say at $0$ and $\infty$, of order two in
$\PSL_2(\IC)$)
we can pull back along the map $w\mapsto z=w^2$
to get a Fuchsian 
system with two apparent singularities at $0$, $\infty$
and four non-apparent singularities at $\pm 1,\pm\sqrt{t}$.
Removing the apparent singularities (using Schlesinger transformations)
yields a 
system $B$ with $A_5$ type $\xi^2\eta^2$, which may be put in the form
\eqref{eq: lin syst} by a coordinate transformation.
Isomonodromic deformations of $A$ correspond to isomonodromic deformations of
$B$, and one can obtain formulae relating the corresponding \PVI\ solutions.
In practice the formulae are much simpler at different (Okamoto
equivalent) values of the parameters (see Ramani et al. \cite{RGT-quad} (2.7)).
We should emphasise that these transformations
are not really symmetries of the family of Painlev\'e VI
equations since the conditions on the parameters restrict us to a 
co-dimension two subset of the four-dimensional parameter space.
Nonetheless they are precisely what is needed to obtain the eight outstanding
icosahedral solutions, since they all have the desired
factor of $a^2$ in their $A_5$ types.
Indeed for these cases, this procedure gives an algebraic relation 
with a solution having half the number of branches;
Examining table 1 of \cite{icosa}
we see solution 31 $\Rightarrow$ solution 44 and in turn
solution 44 $\Rightarrow$ solution 50.
Similarly 
$$
32  \Rightarrow  45  \Rightarrow  51,\qquad
39  \Rightarrow  47,\qquad
40  \Rightarrow  48,\qquad
41  \Rightarrow  49  \Rightarrow  52.$$


The formula of Ramani et al. that we will use 
to construct these outstanding solutions from known solutions
is as follows. (In fact this is the inverse of 
the formula \cite{RGT-quad} (2.7), 
having converted their parameters to our conventions.)

\begin{prop}[\cite{RGT-quad}]
\label{prop: RGT}
Given a solution $(y_0,t_0)$ of \PVI\ with parameters 
of the form 
$\theta=(0,\theta_2,\theta_3,1)$
then, by taking two square roots,
one obtains a new solution $(y,t)$ with parameters
$\theta=(\theta_3,\theta_2,\theta_2,2-\theta_3)/2$
where
$$
y=\frac{(\tau-1)(\eta+1)}{(\tau+1)(\eta-1)},\qquad
t=\left(\frac{\tau-1}{\tau+1}\right)^2
$$
with $\eta^2=y_0, \tau^2=t_0$.
\end{prop}

Note that negating $\tau$ corresponds to the Okamoto symmetry
$(y,t)\mapsto (y/t,1/t)$ and 
negating both $\eta$ and $\tau$  corresponds to
$(y,t)\mapsto (1/y,1/t).$

In practice this will usually lead to an inefficiently parameterised Painlev\'e
curve.
In the cases at hand this may be remedied as follows.
(In the process we will convert the formula to that most directly useful to
us.)
The relation between the Painlev\'e curve $\Pi'$ 
we end up with and the original
curve $\Pi$ may be summarised by the diagram:

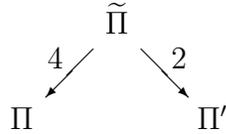
\begin{equation*}
 \setlength{\unitlength}{36pt}
 \begin{picture}(2.2,1.2)(0,0)
 \put(2,0){\object{$\Pi'$}}
 \put(0,0){\object{$\Pi$}}
 \put(1,1){\object{$\wt\Pi$}}
 \put(0.75,0.75){\vector(-1,-1){0.5}}
 \put(0.6,.8){\nwlabel{$4$}}
 \put(1.25,0.75){\vector(1,-1){0.5}}
 \put(1.4,0.8){\swlabel{$2$}}
 \end{picture}
\end{equation*}
where the numbers indicate the degrees of the maps, and $\wt \Pi$  is the
intermediate curve obtained by adjoining the two square roots to the function
field of $\Pi$.

Suppose our initial solution is a pair of functions of the form
\beq \label{eq: initial sol}
Y=\frac{1}{2}+ a_Y(s) u,\qquad T=\frac{1}{2}+a_T(s) u
\eeq
for parameters of the form $\theta=(0,\theta_2,0,\theta_4)$
on a curve of the form
$$\Pi:= \{u^2=u_2(s)\}$$
where $u_2$ is a polynomial, and $a_Y, a_T$ are rational functions of $s$.
In other words $\Pi$ is a double cover of the $s$-line $\IP^1_s$, 
and the symmetry of $\Pi$ 
(negating $u$) corresponds to the symmetry $(y,t)\mapsto (1-y,1-t)$.
Our basic observation is that the parameter $u$ will drop out in the solution
obtained, as follows. 

Applying the Okamoto transformation $(Y,T)\mapsto (Y/(Y-1),T/(T-1))$
yields a solution to which we may apply Proposition \ref{prop: RGT}.
Thus we need to take square roots of $Y/(Y-1), T/(T-1)$, i.e. of expressions
of the form $(A+u)/(A-u)$ where $A=2au_2$ is still a rational function of $s$.
A useful trick is to look for square roots of similar form: i.e. to find 
$B$ such that  
$$\left(\frac{B+u}{B-u}\right)^2=\frac{A+u}{A-u}.$$
Taking the square root of both sides and solving for $B$ we find
$$B=A \pm \sqrt{A^2-u_2}$$
which does not involve $u$.
Carrying this out for both $Y$ and $T$ we obtain 
$$
\eta= \frac{B_Y+u}{B_Y-u},\qquad
\tau= \frac{B_T+u}{B_T-u},
$$ 
where $B_i=A_i \pm \sqrt{A_i^2-u_2}$ for $i=Y,T$.
Then the formulae of Proposition \ref{prop: RGT} yield 
$$y=\frac{B_Y}{B_T},\qquad t=\frac{u_2}{B_T^2}$$
neither of which involves $u$.
Thus $\Pi'$ can be viewed as either the quotient of $\wt \Pi$ by the
involution negating $u$ or as the four-fold cover of the $s$-line obtained by
adjoining functions $v,w$ with 
$$v^2=A_Y^2-u_2,\qquad w^2=A_T^2-u_2$$
where $A_i=2u_2a_i$ for $i=Y,T$.
(The reader may verify that the involution of $\Pi'$ 
negating both $v$ and $w$ together
yields the transformation $(y,t)\mapsto (1/y,1/t)$.)

In turn if we apply the transformation $(y,t)\mapsto (y/(y-1),t/(t-1))$
we will obtain a solution of form similar to \eqref{eq: initial sol}.
In summary (after some relabelling)
the version of the quadratic transformations
we will actually use is as follows:

\begin{cor} \label{cor: qt}
If the functions $y_0,t_0$ of the form
$$ 
y_0=\frac{1}{2}+ a_y(s) u,\qquad t_0=\frac{1}{2}+a_t(s) u
$$
are a \PVI\  solution with parameters $\theta=(0,\theta_2,0,\theta_4)$
on a Painlev\'e curve of the form 
$$\Pi:= \{u^2=u_2(s)\}$$
for a polynomial $u_2(s)$, then the functions
$$
y=\frac{1}{2}+\frac{w+v}{2(A_y-A_t)},\qquad
t=\frac{1}{2}-\frac{A_t}{2w}$$
are a \PVI\  solution for parameters
$\theta=(1-\theta_4,\theta_2,1-\theta_4,2-\theta_2)/2$
on the curve obtained by adjoining to $\IC(s)$ the functions $v,w$ where
$$v^2=A_y^2-u_2,\qquad w^2=A_t^2-u_2$$
and $A_i=2a_iu_2$ for $i=y,t$.
\end{cor}

Of course, a similar result is true upon replacing $\IP^1_s$ by an arbitrary
genus curve, but this will be sufficient for us here.
Note that
 negating both $v$ and $w$ now corresponds to 
the Okamoto transformation $(y,t)\mapsto (1-y,1-t)$.

\end{section}

\begin{section}{Solutions}
We will now carry out the following steps to find the 
formulae for the outstanding 
icosahedral solutions:

1) Choose an icosahedral solution from the 
   table in \cite{icosa} and if possible convert it,
   via Okamoto transformations, into a solution with parameters of the form
   $(0,\theta_2,0,\theta_4)$,

2) Apply Corollary \ref{cor: qt} to obtain a new solution, which will (in
   the examples here) have twice the number of branches (and larger genus) 
   than the original solution,

3) Look for a simple model of the resulting Painlev\'e curve (either as a
   double cover of some $\IP^1$, if it is hyperelliptic, or as a low degree 
   plane curve otherwise). 

A priori 
suitable solutions for step 1) are easily detected by looking for two zero
coordinates in
the solution's alcove point listed in table 1 of \cite{icosa}.

\begin{rmk} \label{rmk: maple file}
To aid the interested reader, and avoid typos, 
a Maple text file of the solutions of this article
has been included with the source file 
(obtained by clicking on ``Other formats'')
for the preprint version on the math arxiv. 
This file also contains solutions equivalent to those written here for which
the corresponding isomonodromic family of 
Fuchsian systems has finite (icosahedral) monodromy group. 
\end{rmk}

\noindent{\bf 10 branch genus zero $\Rightarrow$ 20 branch genus one.}\ \ 

\nobreak

Applying some Okamoto transformations to the $H_3$ solution from 
\cite{Dub95long} E.33, \cite{DubMaz00}, which is equivalent to 
icosahedral solution $32$, one obtains the solution

$$y_0=\frac{1}{2}-{\frac {\left (3\,{s}^{2}+6\,s-1\right )u}{16{s}^{2}}},\quad
t_0=\frac{1}{2}+{\frac {u\, P}{256\,\left (5\,s-1\right ){s}^{3}}}$$
for parameters $\theta=(0,1/5,0,1)$ where $u^2=s$ and
$P=27\,{s}^{5}-315\,{s}^{4}-370\,{s}
^{3}+170\,{s}^{2}-25\,s+1.$
Applying Corollary  \ref{cor: qt} to this (and adjusting $v,w$ slightly to
remove square factors) yields the solution
$$
y=\frac{1}{2}-
{\frac {16\,s \left( 5\,s-1 \right) +vw}
{ 2\left( s-1 \right)  \left( 3\,s+1 \right) v}},\qquad
t=\frac{1}{2}-
{\frac {P}
{ 2\left( s-1 \right)v^2 w}}
$$
for parameters $\theta=(0,1,0,9)/10$ with $P$ as above and where $w=vw_1$ and
\beq \label{eq: soln 45, space curve}
v^2=(9s-1)(s-1),\qquad
w_1^2=s^2-18s+1.
\eeq
One may check directly that this is a genus one solution with twenty branches,
and is equivalent to icosahedral solution $45$.
(It is reassuring to compute the monodromy of the cover $t:\Pi'\to\IP^1$
and find it has the properties listed in table 1 of \cite{icosa}.)
Our next aim is to find a good model of the 
elliptic curve defined by \eqref{eq: soln 45, space curve}, preserving the
symmetry negating $v$.
We will do this by parameterising the conic $w_1^2=s^2-18s+1$ as follows:
$$s= \frac{j^2-1}{2j-18},\qquad w_1=\frac{j^2-18+1}{2j-18}.$$
Then if we define $v=\frac{z}{2j-18}$
the condition that $v^2=(9s-1)(s-1)$
says that $(z,j)$ is a point of the elliptic curve
\beq
z^2= (9j^2-2j+9)(j^2-2j+17),\eeq
and the above formulae give $y,t$ explicitly as functions on this curve.
(One may show, using Magma for example, that this elliptic curve corresponds
to entry 200B1 of Cremona's tables of elliptic curves 
and for example is isomorphic over $\IQ$
to the plane cubic
$u^2= s(s^2 - 5s + 5)$, but this model hides the symmetry of the
Painlev\'e solution.)

Similarly we can proceed with the sibling solution to that above, to obtain 
the solution:
$$
y=\frac{1}{2}-
{\frac {
64\, \left( 5\,s-1 \right) {s}^{2}+\left( s-1 \right)v w
}
{2  \left( 3\,{s}^{3}
+75\,{s}^{2}-15\,s+1 \right)v }}
$$
with $t,s,v,w$ as above but $\theta=(0,3,0,7)/10$. This is equivalent to
icosahedral solution $44$. 

\ 

\noindent{\bf 20 branch genus one $\Rightarrow$ 40 branch genus three.}\ \ 

\nobreak

We can apply Corollary \ref{cor: qt} again 
to the resulting solutions above, since their
parameters are again of the desired form.
Solution $45$ then yields the solution
$$
y=\frac{1}{2}+
{\frac {  \left( {j}^{2}-18\,j+1 \right)z^2 +16\, \left( j+3
 \right)  \left( j+1 \right) vw}
{8 \left( 3\,j-7 \right)  \left( j-9 \right)  \left( j-1 \right) ^{2}v}},
\qquad
t=\frac{1}{2}+{\frac {u\, P}{256\,\left (5\,s-1\right ){s}^{3}}}
$$
with $\theta=(1,1,1,19)/20$,
where $P(s), z^2$ are as in the previous subsection,
$$s=\frac{j^2-1}{2j-18},\qquad u=\frac{w}{2j-18}$$
and now 
\beq\label{eq: 40br space curve}
v^2=-(j-1)(j-9)(5j^2-2j+13),\qquad
w^2=2(j-9)(j^2-1).\eeq

One may check directly that this is a genus 3 solution with forty branches 
and is equivalent to icosahedral solution $51$.
(Note that $t$ is simply the pullback of the original degree $10$
function $t_0$.) 
The curve defined by \eqref{eq: 40br space curve} 
is not hyperelliptic, so we can find
a plane model by taking the canonical embedding.
(Eliminating $s$ from the equations \eqref{eq: 40br space curve}
yields a singular
plane sextic, and we compute three independent differentials directly on this.)
This gives the following model of the Painlev\'e curve as a
smooth plane quartic, with affine equation:
\beq
5(p^4+q^4)+6(p^2q^2+p^2+q^2)+1=0.
\eeq
The solution functions $(y,t)$ become functions on this quartic by setting
$$v = \frac{200p\,(6p^2+5q^2+1)}{84p^2q^2-55q^4-166q^2-156p^2-31},\qquad
w = qv,\qquad  
s = \frac{28v^2-4w^2+800}{3v^2+15w^2-800}.$$
Notice that this curve has three involutions (generating a group isomorphic
to the dihedral group of order eight). These correspond to the
Okamoto symmetries coming from the three hyperplanes on 
which the solution's parameters lie (as
listed in table 1 of \cite{icosa}).
In more detail the symmetries mapping $(p,q)$ to 
$(-p,q),(p,-q),(q,p)$ correspond to the Okamoto symmetries mapping $(y,t)$ to 
$$
(1-y,1-t),\qquad
\left(\frac{y(t-1)}{t-y},1-t\right),\qquad
\left(\frac{y-t}{y-1},t\right)
$$ respectively.

Similarly, from the sibling solution 44 one obtains
the following, which is equivalent to
icosahedral solution $50$:

$$y=\frac{1}{2}+
{\frac { 
\left( {j}^{2}-18\,j+1 \right)  \left( {j}^{2}-2\,j+17 \right) z^{2}
+8\, \left( j-1
 \right)  \left( {j}^{3}+57\,{j}^{2}-69\,j+75 \right) vw}
{8 \left( 3\,{j}^{3}-21\,{j}^{2}-15\,j-31 \right) {w}^{2}v}}
$$
with $t,v,w,z^2$ 
as above and $\theta=(3, 3, 3, 17)/20$.

\ 

\noindent{\bf 15 branch genus one $\Rightarrow$ 30 branch genus two.}\ \ 

\nobreak

If we apply some Okamoto transformations to icosahedral solution 39 (from
\cite{icosa}) then we can obtain the solution
$$ y_0=\frac{1}{2}
-{\frac {u \left( 2\,{s}^{2}+3\,s-3 \right) }
{ 6 \left( s+1\right)  \left( 4\,{s}^{2}+15\,s+15 \right) }},\qquad
t_0=\frac{1}{2}
-{\frac {  u\, P}
{18 \left( 4\,{s}^{2}+15\,s+15 \right) ^{2} \left( {s}^{2}-5 \right) }}
$$
with $\theta=(0, 7/15, 0, 13/15)$, where
$u^2=3(s+5)(4s^2+15s+15)$ and 
$$P= 2\,{s}^{7}+10\,{s}^{6}-90\,{s}^{4}-135\,{s}^
{3}+297\,{s}^{2}+945\,s+675.$$

Applying Corollary \ref{cor: qt} to this, and again adjusting $v,w$ to remove
square factors, yields the solution:
$$
y=\frac{1}{2}+
{
\frac {\left( {s}^{2}-5 \right) {u}^{2}v+ 
s\left( s-3 \right)  \left( s+1 \right) {w}^{3}
}
{ 2\left( s-3 \right)\left( s+5 \right)
\left( {s}^{3}+{s}^{2}-9\,s-15 \right) {w}^{2}}},\quad
t=\frac{1}{2}
+{\frac {   \left( s+5 \right) ^{2}P}
{4s \left( {s}^{2}-9 \right)w^3 }}
$$
with $\theta=(2, 7, 2, 23)/30$ where $P$ and $u^2$ are as above and
\beq \label{eq: g2 space curve}
v^2=s\left( s+5 \right)  \left( s+2 \right)  \left( s-3 \right),\qquad
w^2= s\left( s+5 \right)  \left( s+2 \right)  \left( s+3 \right).
\eeq
This has thirty branches, genus two 
and is equivalent to icosahedral solution $47$.
Being of genus two, the curve \eqref{eq: g2 space curve} is hyperelliptic.
We will express it as a double cover of a $\IP^1$ branched at six points.  
Indeed by choosing a parameter $j$ on the conic $x^2=s^2-9$,
the Painlev\'e curve \eqref{eq: g2 space curve} becomes isomorphic to the 
hyperelliptic curve
$$z^2=(j^2+9)(j+9)(j+1)(j^2+4j+9)$$
via the map
$$v=\frac{j-3}{4j^2}z,\qquad
  w=\frac{j+3}{4j^2}z,\qquad 
  s=\frac{j^2+9}{2j}.$$

Similarly we can repeat starting with solution 40 (the sibling of 39)
and obtain solution 48 (the sibling of solution 47).
The result is 
$$y=\frac{1}{2}
+{\frac { 
\left( {s}^{2}-5 \right) u^2+ 
 \left( {s}^{2}-6\,s-15\right) v w}
{2\,s \left( s+5 \right)\left( s+3 \right) ^{2}{v}}}
$$
with $\theta= (4, 1, 4, 29)/30$ and with $t,v,w,u^2,s$ 
as for solution 47 above.

\ 

\noindent{\bf 18 branch genus one $\Rightarrow$ 36 branch genus three.}\ \ 

\nobreak

Next we will start with icosahedral solution 41 (from
\cite{icosa}; the 10 page implicit form of this solution in the preprint
version of \cite{DubMaz00} is not useful here).
Applying some Okamoto transformations yields the solution:

$$
y_0=\frac{1}{2}-{\frac {8\,{s}^{3}-12\,{s}^{2}+3\,s-4}{6u}},\qquad
t_0=
\frac{1}{2}+
{\frac { P
}
{54\,s\, \left( s-1 \right)u^3 }}.
$$
for $\theta=(0, 1/3, 0, 1)$ where 
\beq\label{eq: 18 br ell curve}
u^2=s(8s^2-11s+8),
\eeq
and
$$P=\left( s+1 \right)  
\left(
32\,({s}^{8}+1)-320\,({s}^{7}+s)+1112\,({s}^{6}+s^2)
-2420\,({s}^{5}+s^3)+3167\,{s}^{4}\right).
$$
Applying Corollary \ref{cor: qt} to this yields the solution:
$$
y=\frac{1}{2}
-{\frac {9\,s \left( s-1 \right) {u}^{2}+ \left( s-2 \right) wv}
{2 \left( {s}^{3}+12\,{s}^{2}-12\,s+4 \right)  \left( 2\,s-1 \right)v}},\quad
t=\frac{1}{2}-
{\frac { P
}
{4 \left( 2\,s-1 \right) \left( s-2 \right)\,v^2\,w}}
$$
with $P,u^2$ as above, $\theta=(0, 1, 0, 5)/6$ and $w=vw_1$ where
\beq \label{eq: g=3 hypell space}
v^2=(s-2)(2s-1)(2s^2+s+2),\qquad w_1^2=s^2-7s+1.\eeq
One may check this is a $36$ branch genus three solution and is equivalent to
icosahedral solution $49$.
However in this case the curve defined by \eqref{eq: g=3 hypell space}
is hyperelliptic.
Indeed let $j$ be a parameter on the conic $w_1^2=s^2-7s+1$, so for example
\beq\label{eq: jparam}
w_1=\frac{j^2-7j+1}{2j-7},\qquad
s= \frac{j^2-1}{2j-7}.
\eeq
Then the Painlev\'e curve \eqref{eq: g=3 hypell space} becomes isomorphic to
$$z^2=\left( {j}^{2}-4\,j+13 \right)\left( 2\,{j}^{2}-2\,j+5 \right) 
\left( 2\,{j}^{4}+2\,{j}^{3}-3\,{j}^{2}-58\,j+107 \right) 
$$
via \eqref{eq: jparam} and the assignment $v=z/(2j-7)^2$.

\ 

\noindent{\bf 36 branch genus three $\Rightarrow$ 72 branch genus seven.}\ \ 

\nobreak

Finally we can apply Corollary \ref{cor: qt} again to the solution above
(since the parameters are of the desired form) to obtain the largest
icosahedral solution.
The solution is given by 
{\tiny
$$y=\frac{1}{2}+
{\frac {9  \left( j-1 \right)  \left( {j}^{3}+27 {j}^{2}-57 j+
79 \right) wv+2\left( 2 {j}^{2}-2 j+5 \right)  \left( {j}^{2}-7 
j+1 \right)  \left( 2 {j}^{4}+2 {j}^{3}-3 {j}^{2}-58 j+107
 \right)  \left( {j}^{2}-4 j+13 \right) ^{2}}{ 6\left( j^2-1 \right)
 \left( 2 {j}^{2}+j+17 \right)  \left( {j}^{3}-3 {j}^{2}+3 j-11
 \right)   \left( 2 j-7 \right) ^{2}v}},
$$}
$$t=
\frac{1}{2}+
{\frac { P}
{54\,s\, \left( s-1 \right)u^3 }},
$$
where $P(s)$ is the polynomial in the previous subsection, 
$$
s= \frac{j^2-1}{2j-7},\qquad
u=\frac{w}{(2j-7)^2},
$$
and
\begin{gather}
v^2=- \left( j+1 \right)  \left( 6+{j}^{2}-2\,j \right)  \left( 4\,{j}^{2}
-13\,j+19 \right),\label{eq: 72 br v2}\\ 
w^2=
\left( j-1 \right)  \left( 2\,j-7 \right)  \left( j+1 \right)
\left( 2\,{j}^{2}+j+17 \right) \left( 4\,{j}^{2}-13\,j+19 \right).
\label{eq: 72 br w2}
\end{gather}
Note that equation \eqref{eq: 72 br w2} is equivalent 
to equation \eqref{eq: 18 br ell curve}
so $t$ is the pullback of the original degree $18$ function $t_0$.

One may check directly that this does indeed define a genus seven, $72$ branch 
Painlev\'e solution and is equivalent to icosahedral solution $52$.
Of course being of genus $7$ the degree--genus formula implies we cannot hope
to find a non-singular plane model of the Painlev\'e curve.
Instead we will look for a low degree plane model with mild singularities.
(The curve obtained upon eliminating 
$j$ from \eqref{eq: 72 br v2}, \eqref{eq: 72 br w2} is a highly singular 
degree $14$ plane curve, with large coefficients.)
We do this by selecting a three-dimensional subspace of the space of 
holomorphic one-forms on the curve, and taking the corresponding plane curve.
After some trial and error choosing a good subspace we found the following 
plane octic with only double points (ten nodes and two tacnodes):
$$
9\,(p^6\,q^2+p^2\,q^6)+
18\,p^4\,q^4+
4\,(p^6+q^6)+
26\,(p^4\,q^2+p^2\,q^4)+
8\,(p^4+q^4)+
57\,p^2\,q^2+
20\,(p^2+q^2)+
16.
$$
The map between the curves is given by
$$
p={\frac {w}{ 3\left( j-1 \right) v}},\qquad
q={\frac {v}{3({j}^{2}-2\,j+6)}}$$  
and, if needed, the (rather long) inverse appears in the accompanying 
computer file (see Remark \ref{rmk: maple file}).

\end{section}

\renewcommand{\baselinestretch}{1}              %
\normalsize
\bibliographystyle{amsplain}    \label{biby}
\bibliography{../thesis/syr}    

\providecommand{\bysame}{\leavevmode\hbox to3em{\hrulefill}\thinspace}
\providecommand{\MR}{\relax\ifhmode\unskip\space\fi MR }
\providecommand{\MRhref}[2]{%
  \href{http://www.ams.org/mathscinet-getitem?mr=#1}{#2}
}
\providecommand{\href}[2]{#2}
\begin{thebibliography}{10}

\bibitem{icosa}
P.~P. Boalch, \emph{The fifty-two icosahedral solutions to {P}ainlev\'e {VI}},
  J. Reine Angew. Math., to appear, (math.AG/0406281, v.7).

\bibitem{srops}
\bysame, \emph{Six results on {P}ainlev\'e {VI}}, math.AG/0503043.

\bibitem{octa}
\bysame, \emph{Some explicit solutions to the {R}iemann--{H}ilbert problem},
  math.DG/0501464.

\bibitem{Dub95long}
B.~Dubrovin, \emph{Geometry of 2{D} topological field theories}, Integrable
  Systems and Quantum Groups (M.Francaviglia and S.Greco, eds.), vol. 1620,
  Springer Lect. Notes Math., 1995, pp.~120--348.

\bibitem{DubMaz00}
B.~Dubrovin and M.~Mazzocco, \emph{Monodromy of certain {P}ainlev\'e-{V}{I}
  transcendents and reflection groups}, Invent. Math. \textbf{141} (2000),
  no.~1, 55--147. \MR{2001j:34114}

\bibitem{Hit-Poncelet}
N.~J. Hitchin, \emph{Poncelet polygons and the {P}ainlev\'e equations},
  Geometry and analysis (Bombay, 1992), Tata Inst. Fund. Res., Bombay, 1995,
  pp.~151--185. \MR{97d:32042}

\bibitem{Kitaev-dessins}
A.~V. Kitaev, \emph{Dessins d'enfants, their deformations and algebraic the
  sixth {P}ainlev\'e and {G}auss hypergeometric functions}, nlin.SI/0309078,
  v.3.

\bibitem{Kitaev-quad-p6}
\bysame, \emph{Quadratic transformations for the sixth {P}ainlev\'e equation},
  Lett. Math. Phys. \textbf{21} (1991), no.~2, 105--111.

\bibitem{OkaPVI}
K.~Okamoto, \emph{Studies on the {P}ainlev\'e equations. {I}. {S}ixth
  {P}ainlev\'e equation {$P\sb {{\rm VI}}$}}, Ann. Mat. Pura Appl. (4)
  \textbf{146} (1987), 337--381. \MR{88m:58062}

\bibitem{RGT-quad}
A.~Ramani, B.~Grammaticos, and T.~Tamizhmani, \emph{Quadratic relations in
  continuous and discrete {P}ainlev\'e equations}, J. Phys. A \textbf{33}
  (2000), no.~15, 3033--3044. \MR{MR1766506 (2001d:34018)}

\bibitem{TOS-folding}
T.~Tsuda, K.~Okamoto, and H.~Sakai, \emph{Folding transformations of the
  {P}ainlev\'e equations}, Math. Ann. \textbf{331} (2005), 713--738.

\end{thebibliography}
\end{document}